\documentclass[conference,a4paper]{IEEEtran}
\usepackage[pdftex]{graphicx}
\graphicspath{{../pdf/}{../jpeg/}}
\DeclareGraphicsExtensions{.pdf,.jpeg,.png}

\usepackage{amsmath}
\usepackage{multirow}
\usepackage{resizegather}

\usepackage{array}
\usepackage{mdwmath}
\usepackage{mdwtab}
\usepackage{eqparbox}
\usepackage{mathtools}

\usepackage{url}

\usepackage[table,xcdraw]{xcolor}
\usepackage{placeins}

\usepackage[utf8]{inputenc} % allow utf-8 input
\usepackage[T1]{fontenc}    % use 8-bit T1 fonts
\usepackage{url}            % simple URL typesetting
\usepackage{booktabs}       % professional-quality tables
\usepackage{amsfonts}       % blackboard math symbols
\usepackage{amsmath}       % blackboard math symbols
\usepackage{nicefrac}       % compact symbols for 1/2, etc.
\usepackage{microtype}      % microtypography
\usepackage{lipsum}
\usepackage{fancyhdr}       % header
\usepackage{graphicx}       % graphics
\usepackage{amssymb}

\graphicspath{{media/}}     % organize your images and other figures under media/ folder

\usepackage{algorithm}
\usepackage{algpseudocode}
\usepackage{xurl} % not \usepackage{url}
\usepackage[hidelinks]{hyperref}

\DeclareMathOperator*{\argmin}{\text{arg}\,min}

\usepackage{color}
\usepackage{tabularray}
\definecolor{PeriwinkleGray}{rgb}{0.749,0.839,0.909}
\definecolor{PastelGreen}{rgb}{0.419,0.901,0.603}
\definecolor{PetiteOrchid}{rgb}{0.85,0.588,0.588}

\begin{document}
\title{Controlling Large Electric Vehicle Charging Stations via User Behavior Modeling and Stochastic Programming}

\author{\IEEEauthorblockN{Alban Puech$^{1,2}$, Tristan Rigaut$^1$, William Templier$^1$, Maud Tournoud$^1$} 
\IEEEauthorblockA{\textit{$^1$ Schneider Digital - AI Hub, Schneider Electric, Grenoble, France} \\
\textit{$^2$ École Polytechnique Fédérale de Lausanne (EPFL), Lausanne, Switzerland}\\
alban.puech@epfl.ch}
}

\pagestyle{fancy}
\thispagestyle{empty}
\rhead{ \textit{ }} 

\linespread{0.9}

% Update your Headers here
\fancyhead[LO]{Controlling Large Electric Vehicle Charging Stations via User Behavior Modeling and Stochastic Programming}
% \fancyhead[RE]{Firstauthor and Secondauthor} % Firstauthor et al. if more than 2 - must use \documentclass[twoside]{article}

% ====================================================================

\IEEEoverridecommandlockouts
\IEEEpubid{\makebox[\columnwidth]{979-8-3503-9042-1/24/\$31.00~\copyright2024 IEEE \hfill}
\hspace{\columnsep}\makebox[\columnwidth]{ }}
\maketitle

% === ABSTRACT ====================================================================
% =================================================================================
\begin{abstract}
This paper introduces an Electric Vehicle Charging Station (EVCS) model that incorporates real-world constraints, such as slot power limitations, contract threshold overruns penalties, or early disconnections of electric vehicles (EVs). We propose a formulation of the problem of EVCS control under uncertainty, and implement two Multi-Stage Stochastic Programming approaches that leverage user-provided information, namely, Model Predictive Control and Two-Stage Stochastic Programming. The model addresses uncertainties in charging session start and end times, as well as in energy demand. A user's behavior model based on a sojourn-time-dependent stochastic process enhances cost reduction while maintaining customer satisfaction. The benefits of the two proposed methods are showcased against two baselines over a 22-day simulation using a real-world dataset. The two-stage approach demonstrates robustness against early disconnections by considering a wider range of uncertainty scenarios for optimization. The algorithm prioritizing user satisfaction over electricity cost achieves a 20\% and 36\% improvement in two user satisfaction metrics compared to an industry-standard baseline. Additionally, the algorithm striking the best balance between cost and user satisfaction exhibits a mere 3\% relative cost increase compared to the theoretically optimal baseline --- for which the nonanticipativity constraint is relaxed --- while attaining 94\% and 84\% of the user satisfaction performance in the two used satisfaction metrics.

\end{abstract}

% \fi
\IEEEpeerreviewmaketitle

% ====================================================================
% ====================================================================
% ====================================================================

% === KEYWORDS ====================================================================
% =================================================================================

\section{Introduction}
% \blfootnote{$^\dagger$ Work carried out while at Institut Polytechnique de Paris, Palaiseau, France}

The number of public electric vehicle (EV) charging stations (EVCS) increased by 55\% in 2022~\cite{IEA}, and this number is expected to grow as the EU and the US multiply incentives to support the development of publicly accessible fast-charging infrastructure~\cite{USplans, EUplans}. However, EV charging stations come with many challenges, both in their operation (e.g. load management, user access, payment, or compatibility) and in their integration into the existing electrical grids~\cite{challenges}. Additionally, EVCSs typically have to satisfy power constraints, while satisfying the requests of their users and minimizing electricity costs. In this context, smart control strategies are required and will play a key role in scaling the EVCS network.

The difficulty of EVCS control comes from the uncertainties involved in their operation. While users can be asked to indicate their intended parking time, they often, in practice, disconnect their vehicles much before; creating a need for session-end forecasts. As EVCSs frequently run on contracts having penalties for consumption peaks, planning the charge requires modeling the arrival time, session duration, and energy requests of future sessions. This turns the problem of EVCS control into a stochastic optimization problem, which can be tackled using different Stochastic Optimization approaches.
In \cite{9960988}, Cheng et al. proposed an EV charging algorithm to minimize carbon emissions while meeting user demands. Unlike their approach, we predict arrival and departure times to handle early disconnections and better model EVCS load. We also consider constraints on maximum power delivery, addressing real-world conditions \cite{POWELL2020115352, muratori2018impact}. Tucker et al. in \cite{9917194} model transformer capacity and forecast EV departures but determine departure times deterministically. In contrast, our method models connection and departure times based on sojourn times, incorporating features like charging slot indices. We also ensure customer satisfaction through objective function embedding, even with early disconnections. The necessity for charging session forecasting is discussed in \cite{GENOV2024121969}, advocating for machine learning-based prediction, including weather and traffic data. Our approach advances charging behavior modeling by generating connection and disconnection scenarios for stochastic optimization. More importantly, we do so by leveraging locally available data, while fully aligning with implementation constraints \cite{tristanspaper}.

In this paper, we thus propose an approach that minimizes charging costs while preserving user satisfaction. The solution can easily be implemented on an industrial system and leverages the requests of the users while anticipating early disconnections. Our contributions are:
\begin{itemize}
    \item  An EVCS model that incorporates real-world constraints, such as slot power limitation, contract threshold overruns penalties, or early disconnection of EVs.
    \item A Two-Stage Stochastic Programming approach that utilizes charging session scenarios instead of a single forecast of future charging sessions.
    \item A modeling of the charging sessions' end and start, formulated as a sojourn-dependent stochastic process, which improves the modeling approach presented in~\cite{tristanspaper} and leverages user-provided information.
    \item A comparative analysis of our novel control method against two baselines in a simulation on a 22-day real-world dataset of an EVCS, comprising 32 slots.
\end{itemize}

In section \ref{sec:statement}, we introduce the EVCS model with its dynamics and parameters. In section \ref{sec:opti}, we state the EVCS control problem as a generic stochastic optimization problem. We describe the Multi-Stage-Stochastic Programming algorithms, along with our user behavior modeling, in \ref{sec:methodology}. Finally, we present and discuss the numerical results obtained through simulations in section \ref{sec:results}.

\section{EVCS Model}
\label{sec:statement}

We model the EVCS as a discrete dynamical system. Control actions are taken at each time step $t$, of length $\Delta t$, from the sequence $\mathcal{T} = \{0,1,\ldots,T\}$. We denote by $n \in \mathbb{N}^*$ the number of slots. Slots are independent charging points that users connect their EVs to. The binary variable ${o^i_t} \in \{0,1\}$ indicates whether the slot $i,~i<n$ is in an ``active'' or ``inactive'' state at time $t\in \mathcal{T}$.
\subsection{Charging session start} 
When a customer connects their EV to an ``inactive'' slot, they hand over a request to the slot controller, as per~\cite{iso}. A request is a pair $(k, \delta)$, where $k \in \mathbb{R}_{\geq0}$ is the amount of electrical energy (kWh) to be provided to the EV --- net of any transmission or charging loss --- and $\delta \in \mathbb{N}^*$ is the announced parking time, expressed as a number of time steps. In practice, one can assume $k$ to be directly given by the user, or obtained through an EV-companion app that either has access to the current state of charge of the EV or converts a mileage capacity request into an energy request. This connection starts a charging session, and the associated slot state switches from ``inactive'' to ``active''. The binary variable  ${a^i_t}  \in \{0,1\}$ takes a positive value if slot $i$ gets connected to an EV, receives a request, and thus switches from ``inactive'' to ``active'' state at time step $t$. This variable is activated only at connection time.

\subsection{Charging session end} 
An ``active'' charging slot remains in its state until the charging session ends. A session end is triggered if the EV gets disconnected, \textbf{or if the announced parking time elapses}, even if the EV stays connected. Users may leave their EV connected to the slot longer than the announced end of parking time --- resulting in idle time. The session end is thus designed in such a way that the controller cannot rely on this idle time to satisfy the user request, so as to discourage users from keeping their EV parked longer than initially indicated. %At time $t$, the binary variable ${d^i_t\in \{0,1\}}$ takes a positive value if the EV connected to slot $i$ gets disconnected by its owner. We similarly use 
We use the binary variable ${q^i_t\in \{0,1\}}$ to indicate if a slot switches from an ``active'' to an ``inactive'' state at time $t$. At the controller level, there is no distinction between ``inactive'' slots that are still connected to an EV (e.g., after the announced parking time elapsed), and free slots.

\subsection{Charging dynamics}
At each time step $t\in \mathcal{T}$, an ``active'' slot $i,~i<n$ can deliver an energy amount $e^i_t \leq \bar{e},~e^i_t \in \mathbb{R}_{\geq0}$, where $\bar{e}\in \mathbb{R}_{>0}$ denotes the maximum energy amount that can be delivered by a slot in the span of a single time step. The remaining energy $r^i_{t} \in \mathbb{R}_{\geq0}$ to provide the EV connected to the ``active'' slot $i$ with, at time step $t$, is given by:
\begin{equation}
r^i_{t} = \left\{ \,
\begin{IEEEeqnarraybox}[][c]{l?s}
\IEEEstrut
    k^i_t \text{ if } t=t_0 \\
    r^i_{t-1} - \eta \cdot e^i_{t-1} \text{ if } t_0 < t \leq t_0 + \delta^i_{t_0}
\IEEEstrut
\end{IEEEeqnarraybox}
\right.
\end{equation}

where $\eta \in (0,1]$ denotes the charge efficiency and $t_0 \in \mathcal{T}$ the time step during which the last request $(k^i_{t_0}, \delta^i_{t_0})$ was made.

At each time step $t \in \mathcal{T}$, for each ``active'' slot $i$, considering a session with request $(k^i_{t_0}, \delta^i_{t_0})$ starting at $t_0$, we also introduce a variable $z^i_t$ which stores the initial energy request of the session, that is, $z^i_{t} = k^i_{t_0}~\forall t,~t_0 \leq t \leq t_0 + \delta^i_{t_0}$.
This variable will prove useful later to measure the satisfaction of the users. 

Finally, we store the remaining number of time steps from the current time step $t$ before the announced end-of-parking time of slot $i$. We call this feature $m_t^i\in (\mathbb{N} \cup \{-\infty\})$.

\subsection{Costs and objectives}
We control the EVCS in order to reach a trade-off between customer satisfaction and costs.

\subsubsection{Customer satisfaction} We define the satisfaction score $\chi^i_t$ of an ``active'' slot $i,~i<n$ at time $t$ using the ratio between the remaining energy $r_i^t$ to provide to the EV, and the associated initial energy requested $z^i_{t}=k^i_{t_0}$, where $t_0$ is the time step corresponding to the last connection:
\begin{equation}\label{eq:satisfaction}
\chi^i_t = 1-\frac{r_i^t}{z^i_{t}}
\end{equation}

\subsubsection{Electricity costs} The EVCS consumption at time $t \in \mathcal{T}$ is denoted by $c_t \in \mathbb{R}_{\geq0}$ and corresponds to the sum of the consumption of its slots. The electricity cost at time step $t$ is denoted by $\zeta_t\in \mathbb{R}_{\geq0}$ and is expressed as follows, with $p_t \in \mathbb{R}_{\geq0}$ denoting the electricity price at time $t$:
\begin{equation}\label{eq:elec_cost}
\zeta_{t} = p_t \times c_t, \qquad c_t = \sum_{i\in S} e^i_{t}
\end{equation} 

\subsubsection{Grid import threshold overrun} At every time step $t \in \mathcal{T}$, a fixed penalty $\xi \in \mathbb{R}_{>0}$ is paid if the EVCS electricity consumption exceeds a threshold $\bar{c} \in \mathbb{R}_{>0}$, that is, if  $c_t > \bar{c}$. At each time step $t$, the cost incurred by this threshold is expressed as follows, where $\mathbf{1}_{> \bar{c}}: \mathbb{R} \to \{0,1\}$ is the indicator function that takes the value $1$ when its argument is above $\bar{c}$ and $0$ otherwise:
\begin{align} \label{eq:penalties}
\mathbf{1}_{> \bar{c}}(c_t) \times \xi
\end{align}

\section{Optimal control formulation}
\label{sec:opti}
\subsection{State, exogenous uncertainty and control variables} We formulate hereunder the problem as a generic Stochastic Optimal Control problem using standard notations used in the Control literature~\cite{bertsekas2012dynamic}. The \textbf{state variable} at time $t \in \mathcal{T}$ is denoted by $x_t\in \mathbb{X} = \big(\{0,1\} \times \mathbb{R}_{\geq0} \times \mathbb{R}_{\geq0} \times (\mathbb{N} \cup \{-\infty\})\big)^n$: 
\begin{align}
x_t &= (\underbrace{o^i_t}_{\substack{\text{state} \\ \text{$\in$} \\ \text{"active",} \\ \text{"inactive"}}}, \underbrace{r^i_t}_{\substack{\text{remaining} \\ \text{energy} \\ \text{to provide}}}, \underbrace{z^i_t}_{\substack{\text{initial energy} \\ \text{request of} \\ \text{the session}}}, \underbrace{m^i_t}_{\substack{\text{remaining number} \\ \text{of steps before} \\ \text{announced end} \\ \text{of parking}}})_{i=0}^{n-1}
\end{align}
The \textbf{exogenous uncertainty variable} is denoted by $w_t~\in~\mathbb{W}$ where $\mathbb{W}=\big(\{0,1\}^2 \times~\mathbb{R}_{\geq0}~\times~\mathbb{N^*}~\big)^n$:
\begin{align}
    w_t &= (\underbrace{a^i_t}_{\substack{\text{session} \\ \text{start} \\ \text{boolean}}}, \underbrace{q^i_t}_{\substack{\text{session} \\ \text{end} \ \text{boolean}}}, \underbrace{k^i_t}_{\substack{\text{initial} \\ \text{kwh} \ \text{request}}}, \underbrace{\delta^i_t}_{\substack{\text{announced} \\ \text{parking} \\ \text{time}}})_{i=0}^{n-1}
\end{align}
The \textbf{control (action)} at time $t \in \mathcal{T}$ is denoted by $u_t = (e^i_t)_{i=0}^{n-1}$ in the set $\mathbb{U}=(\mathbb{R}_{\geq0})^n$ and gathers the power provided by each slot. The \textbf{dynamics} and the \textbf{constraints} are expressed as $x_{t+1} = f_t(x_t, u_t, w_t)$ and $g_t(x_t,u_t,w_t) \leq 0$, where $f_t: (\mathbb{X},\mathbb{U},\mathbb{W}) \to \mathbb{X}$ and $g_t:(\mathbb{X},\mathbb{U},\mathbb{W}) \to \mathbb{R}^{n_c}$ denote respectively the dynamics function and the constraints, and where $n_c$ denotes the number of constraints. We define the \textbf{stage cost} $L_t: (\mathbb{X},\mathbb{U}) \to \mathbb{R}_{\geq0}$ as a linear combination of the electricity cost~\eqref{eq:elec_cost}, the grid threshold overrun penalty~\eqref{eq:penalties} and the satisfaction score~\eqref{eq:satisfaction}, with weight given to the satisfaction denoted $\alpha \in \mathbb{R}_{>0}$. This weight gives a monetary value to the satisfaction of the users:
\begin{equation}
L_t(x,u) = \underbrace{\big( \zeta_{t} + \mathbf{1}_{>\bar{c}}(c_t) \times \xi\big)}_{\text{electricity cost including penalties}} + ~\alpha \cdot \sum_{i<n}\underbrace{(1-\chi^i_t)}_{\substack{\text{customer} \\ \text{dissatisfaction}}}
\end{equation}

In general, this stage cost can depend on the exogenous uncertainties. We thus assume that $L_t$ is a function mapping $(\mathbb{X},\mathbb{U}, \mathbb{W})$ to $\mathbb{R}_{\geq0}$.

\subsection{Controller definition and Multi-Stage Stochastic Programming policies}
Our goal is to define a control strategy or policy. A  policy is a collection 
$\Phi = \{\phi_{t_0}\}_{{t_0} \in \mathcal{T} }$
of mappings (one for each ${t_0} \in \mathcal{T}$), that each takes a current state $x_{t_0} \in \mathbb{X}$, and an uncertainty history $h_{t_0} =\{w_t\}_{t \leq {t_0}}  $, and that returns a control action $u_{t_0}$ solution to the following Multi-Stage Stochastic Programming problem:
\begin{subequations}
\begin{align}
    \argmin_{u_{t_0} \in \mathbb{U}}&  \min_{(\mathbf{U}_{t_0+i})_{1\leq i\leq R}} \mathbb{E} \sum_{t=t_0}^{t_0+R} L_{t}( \mathbf{X}_{t}, \mathbf{U}_{t}, \mathbf{W}_{t})  \label{eq:expectation} \\
    \text{s.t. }& \mathbf{X}_{t+1} = f_{t}(\mathbf{X}_{t}, \mathbf{U}_{t}, \mathbf{W}_{t}) \\
    & g_{t}(\mathbf{X}_{t}, \mathbf{U}_{t}, \mathbf{W}_{t}) \leq 0\\
    & \mathbf{X}_{t_0}=x_{t_0}, \mathbf{H}_{t_0}=h_{t_0}\\
    & \sigma(\mathbf{U}_t) \subset \sigma(\mathbf{H}_{t_0},\mathbf{W}_{t_0+1},\ldots,\mathbf{W}_{t})
\end{align}
\end{subequations}
where $R\in \mathbb{N}^*$ is a control horizon and where the expectation is computed over the conditional distribution of the uncertainties given the uncertainty history $\mathbb{P}(\mathbf{W}_{t_0}, \mathbf{W}_{t_0+1},\ldots,\mathbf{W}_{t_0+R} | \mathbf{H}_{t_0} = h_{t_0})$. 

\section{Methodology}
\label{sec:methodology}
As the uncertainty distribution has infinite support, we solve the previously formulated problem on a limited number of uncertainty scenarios, obtained from a scenario tree. We thus have a tractable approximation of the expectation displayed in \eqref{eq:expectation}. We here refer to~\cite{puech2023controlling} for further details about scenario trees, that conceptualize the strategy adopted here.

\subsection{Multi-Stage Stochastic Programming algorithms}

We use two Multi-Stage Stochastic Programming (MSSP) algorithms, namely, Two-Stage Stochastic Programming (\texttt{2S}) and Model Predictive Control (\texttt{MPC}).

\subsubsection{Two-stage Stochastic Programming (\texttt{2S})}

Two-stage Stochastic Programming uses a collection of~$K\in \mathbb{N^*}$ samples indexed by~$0 \leq k < K$ denoted ($w_{t_0},\tilde{w}^k_{t_0+1},...,\tilde{w}^k_{t_0+R})_{k=0}^{K-1}$ drawn from the uncertainty distribution $\mathbb{P}(\mathbf{W}_{t_0},\ldots,\mathbf{W}_{t_0+R} | \mathbf{H}_{t_0} = h_{t_0})$, where the first variable of the sampled tuples, $w_{t_0}$, is observable and given by the history $h_{t_0}$ and is thus fixed. These samples are then assigned to $K' \leq K$ clusters. The probability $\pi_k \in (0,1]$ of a cluster is then computed as the ratio between the number of samples associated with this cluster, and the total number of samples. For each of the $K'$ cluster centers, we find their closest sample, and we denote by ($w_{t_0},\hat{w}^k_{t_0+1},...,\hat{w}^k_{t_0+R})_{k=0}^{K'-1}$ this new collection, that we then use to approximate the previously formulated problem by the following smaller and tractable optimization problem:
\begin{subequations} 
\begin{align*}
    \argmin_{u_{t_0}\in U}  &\min_{{(u^k_t)_{\substack{t_0 < t \leq t_0+R \\ k < K'}}}}  \sum_{\substack{t_0\leq t \leq t_0+R \\  k < K'}}  \pi_k \cdot L_{t}(x^k_{t}, u^k_{t}, \hat{w}^k_{t}) \\
    \text{s.t. }& x^k_{t+1} = f_{t}(x^k_{t}, u^k_{t}, \hat{w}^k_{t}),~\forall k < K',~\forall t,~t_0\leq t < t_0+R \\
    &g_{t}(x^k_{t}, u^k_{t}, \hat{w}^k_{t})\leq 0,~\forall k < K',\quad ~\forall t,~t_0\leq t \leq t_0+R 
    % & x^k_{t_0} = x_{t_0},\quad \forall k < K'
\end{align*}
\end{subequations}

\subsubsection{Model Predictive Control (\texttt{MPC})}
\label{sec:MPC}
Model Predictive Control uses a unique uncertainty forecast $({w}_{t_0},\hat{w}_{t_0+1},...,\hat{w}_{t_0+R})$ of the future realizations of $(\mathbf{W}_{t_0},\ldots,\mathbf{W}_{t_0+R})$ conditioned by $\mathbf{H}_{t_0} = h_{t_0}$. In practice, the control action at time $t_0$ is thus obtained by solving the following optimization problem:
\begin{subequations} 
\begin{align*}
    \argmin_{u_{t_0}\in U}  &\min_{{(u_t)}_{t_0 < t \leq t_0+R}}  \sum_{t_0 \leq t \leq t_0+R}  L_{t}(x_{t}, u_{t}, \hat{w}_{t}) \\
    \text{s.t. } & x_{t+1} = f_{t}(x_{t}, u_{t}, \hat{w}_{t}),~ \forall t,~t_0\leq t < t_0+R \\
    &g_{t}(x_{t}, u_{t}, \hat{w}_{t})\leq 0,\quad ~\forall t,~t_0\leq t \leq t_0+R
\end{align*}
\end{subequations}

\subsection{Uncertainty distribution modeling}
In the following subsection, we explain the uncertainty distribution modeling necessary to obtain the uncertainty samples used with \texttt{2S}, and the forecasts used with \texttt{MPC}.
%TODO finish reading here
\subsubsection{Session start and end}
We model the charging session state transition of each slot $i,~i<n$ by a sojourn-dependent stochastic process $\mathbf{O}^i$ with two possible states, $\mathbf{O}^i(t) = 0$  (inactive) or $\mathbf{O}^i(t) = 1$ (active). This formulation is inspired by the Sojourn-time-dependent semi-markov switching processes introduced in~\cite{campo}. Sojourn-time-dependent semi-markov switching processes are semi-markov processes in which the switching (i.e. transition between two distinct states) probability depends on the sojourn time (i.e. the time elapsed since the last switch). In our approach, we additionally condition the switching probability on the following set of features, which we directly derive from the uncertainty history $h_t$: 
\begin{itemize}
    \item \textit{Sojourn time.}
We store the number of time steps since the last transition from ``inactive'' to ``active'' states, or from ``inactive'' to ``active'' states, in the variable $g_t^i \in \mathbb{N}$.
    \item \textit{Time related-features.} An encoding of the current hour of the day, and the weekday, denoted by $l_t \in \mathbb{N}$. %TODO: specify what set
    \item \textit{Slot index.} Denoted $\tilde{i},~\tilde{i}<n$.
\end{itemize}

The transition probability $p_{y\to z}$ from state $y \in \{0,1\}$ to state $z \in \{0,1\}$ of slot $i$ is thus defined as: $p^i_{y\to z} = \mathbb{P}\big(\mathbf{O}(t) = z \big| \mathbf{O}(t-1) = y, \underbrace{g=g^i_{t-1}}_{\text{sojourn time}}, l=l_{t-1}, \tilde{i}=i\big)$

\subsubsection{kWh request}
\label{sub:kwh}
% Likewise for a given slot $i$, sojourn time $g$, remaining time before announced departure $m$ and time of the day $t$ we could model the distribution of possible kWh requests. Instead, we decided to ignore the stochasticity of the request once the sojourn time, the remaining time, the time of the day, and the index of the slot are known.

We model the kWh request $k^i_t$ as a deterministic function $\mathbf{k}$ of the sojourn time $g$, the time of the day $l$, and the slot index $i$: $\mathbf{k}: (\mathbb{N} \times \mathbb{N} \times \{0,1,...,n-1\}) \to \mathbb{R}_{\geq0} $. %TODO: correct sets

\subsection{Uncertainty distribution approximation}
\label{sub:uncertainty}

To approximate the uncertainty distribution and generate uncertainty samples ($w_{t_0},\tilde{w}^k_{t_0+1},...,\tilde{w}^k_{t_0+R})_{k=0}^{K-1}$, we use a training dataset containing past charging sessions to fit two gradient-boosting classifiers and a regressor. These are then used to obtain each of the uncertainty components:
    \subsubsection{Session starts} The first classifier $c_1$ is used to obtain the transition probabilities from ``inactive'' to ``active'' state, so as to obtain the session start uncertainty samples denoted $(a_{t_0},\tilde{a}_{t_0+1},...,\tilde{a}_{t_0+R})$. It approximates the distribution $p^i_{0\to 1} = \mathbb{P}(\mathbf{O}(t) = 1 | \mathbf{O}(t-1) = 0, g=g^i_{t-1}, l=l_{t-1}, \tilde{i}=i)$ conditioned by the previous slot activation state, the sojourn time, the time and date, and the slot index.
    \subsubsection{Kwh requests} The regressor $r_1$ is used to obtain the kWh request associated with a session start (transition from ``inactive'' to ``active'' state). It approximates $\mathbf{k}$, the function that returns the kWh requests given $(g, l, \tilde{i})$. 
    \subsubsection{Session ends} The second classifier is used to obtain the transition probabilities from the ``active'' to ``inactive'' state, and thus, the session end uncertainty samples \sloppy $(q_{t_0},\tilde{q}_{t_0+1},...,\tilde{q}_{t_0+R})$. It approximates: $p^i_{1\to 0 } = \mathbb{P}(\mathbf{O}(t) = 0 | \mathbf{O}(t-1) = 1, g=g^i_{t-1}, l=l_{t-1}, \tilde{i}=i)$
    % \item The second and third classifiers are used to obtain the transition probabilities from the ``active'' to ``inactive'' state, and thus, the session end uncertainty samples $(\tilde{q}_{t_0},\tilde{q}_{t_0+1},...,\tilde{q}_{t_0+R})$. 
    % \begin{itemize}
    % \item The second one is used on ``active'' slots for which user requests are available, that is, $m > -\infty$ (typically all charging sessions, except those obtained using $c_1$). It is denoted $c_2$ and approximates  :   $$p^i_{1\to 0} = \mathbb{P}\bigg(O(t) = 1 \bigg| O(t-1) = 0, g=g^i_{t-1},m=m^i_{t-1}, l=l_{t-1}, \tilde{i}=i\bigg)$$
    % \item The third one is used to predict session end of sessions for which $m=-\infty$ (future charging sessions predicted using $c_1$. It is denoted $c_3$ and approximates:   $$p^i_{1\to 0 } = \mathbb{P}\bigg(O(t) = 1 \bigg| O(t-1) = 0, g=g^i_{t-1}, l=l_{t-1}, \tilde{i}=i\bigg)$$
    % \end{itemize}

\begin{figure}[htp]
    \centering
\includegraphics[height=0.3\textheight]{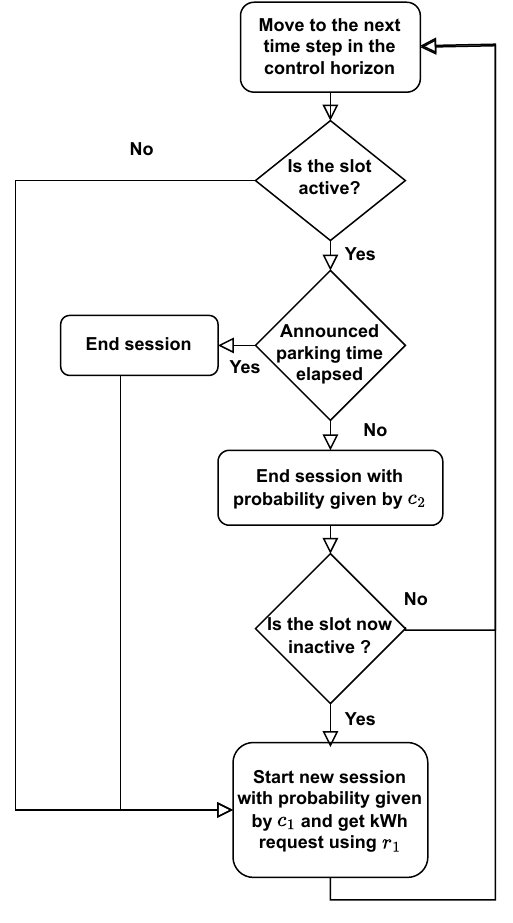}
    \caption{Working principle of a single iteration of the scenario generation algorithm on a single charging slot.}
    \label{fig:flow_chart}
\end{figure}

The working principle of the algorithm used to generate the uncertainty samples using $c_1, c_2$ and $r_1$ is illustrated in Fig. \ref{fig:flow_chart}
%\footnote{The complete pseudo-code of the algorithm will be provided in the appendix of the camera ready version}. This modeling shows numerous improvements over the one introduced in~\cite{tristanspaper}
. Regarding the end-of-session times, our novel model allows for correcting the tendency of users to disconnect their vehicles before the announced end of parking time, while~\cite{tristanspaper} was assuming this announced time to be the effective end-of-session time. In~\cite{tristanspaper}, the controller could plan on charging the EVs on the very last time steps of the announced parking period, leading to customer dissatisfaction if the users disconnected their EVs earlier. By allowing our model to predict session ends to happen before the expected time, this risk is reduced. Moreover, as opposed to the heuristic proposed in \cite{tristanspaper} consisting in setting the load of free slots to constant values to simulate future charging sessions; planning the sessions start, sessions end, and kWh request allows for more accurate modeling of the future EVCS load.

\subsection{Uncertainty forecasts for \texttt{MPC}}
To obtain different uncertainty samples (used with \texttt{2S}), the switches between states are done with probability obtained using $c_1$ and $c_2$, allowing different scenarios to be obtained for the same state conditions. %As introduced in \ref{sec:MPC}, our \texttt{MPC} algorithm uses an uncertainty forecast instead of relying on a collection of uncertainty samples. In this context, the switch between states is done in a deterministic way, i.e. only if the probability of state switch (obtained using $c_1$ and $c_2$) is greater than the probability of staying in the same state. In practice, this means that lines 13, and 25 are modified to compare p with 0.5 instead of a random value.
     
% It means that, considering a slot $i$ at each step $t$, the algorithm takes as input the current state of the slot: active/inactive, sojourn time and time of the day. If the slot is active, we take the probability of transitioning to an inactive state from the classifier $c_1$. if this probability is above $0.5$, the transition is made and the algorithm continues to explore the next time steps to build the most probable sequence of transitions.

%TODO:Started reading again

\subsection{Baseline algorithms}

We compare our two MSSP algorithms (\texttt{MPC} and \texttt{2S}) against two baselines.

\subsubsection{User request-based MPC (\texttt{R-MPC})}
The first one, which we call \texttt{R-MPC} is a Model Predictive Control approach that solely relies on the user requests to build the uncertainty forecast. It is the algorithm introduced in~\cite{tristanspaper}. The user-announced end-of-parking time \textbf{is assumed to be true} and \textbf{no modeling of the future charging session start and end is done}. Instead, the load on empty slots is replaced by an uncontrollable average energy consumption given the time of the day and the slot identifier. These average energy consumption values are obtained by running \texttt{P-MPC} (as defined in the next section) on a training set. We recall the expression of the exogenous uncertainties:
\begin{equation}
 w_t =  (\underbrace{a^i_t}_{\text{session start}}, \underbrace{q^i_t}_{\text{session end}}, \underbrace{k^i_t}_{\text{ initial kWh request}}, \underbrace{\delta^i_t}_{\substack{\text{announced} \\ \text{parking} \\ \text{time}}})_{i=0}^{n-1},
\end{equation}
The session ends of currently ``active'' slots are computed using the announced parking times. Moreover, no future session is predicted, so that $a^i_t = k^i_t = \delta^i_t = 0 ~\forall t,~t_o+1 \leq t \leq t_o + R$. 

% Assuming that slot $i$ is ``active'' and that the last request $r_{t_0}^i = (k_{t_0}^i, \delta_{t_0}^i)$, we get: 
% \begin{equation}
% q^i_{t} =  \left\{ \,
% \begin{IEEEeqnarraybox}[][c]{l?s}
% \IEEEstrut
%     1 \text{ if } t=t_0+\delta_{t_0}^i,\\
%    0 \text{ if } \neq t_0+\delta_{t_0}^i
% \IEEEstrut
% \end{IEEEeqnarraybox}
% \right.
% \end{equation}

 %This mean load corresponds to the one under the \texttt{P-MPC} controller ran on a training dataset containing data for the same slots but on a different period than the one covered by the test set.

This baseline is the most straightforward MPC approach that can be deployed using the available user data. However, it does not anticipate the end of charging sessions by disconnections before the announced parking time elapses and does not account for potential future charging sessions that could increase the load in future time steps. This becomes problematic if the EVCS delays the charge of the current EVs on those time steps --- e.g., because of off-peak prices --- as it results in unsatisfied charging requests.

\subsubsection{\texttt{P-MPC}: Perfect MPC}

The second baseline, called \texttt{P-MPC}, is an MPC algorithm for which the nonanticipativity constraint is relaxed. In practice, this means that it is given the true realizations of the random variables as forecast.

\section{Results}
\label{sec:results}

We present numerical results obtained by implementing the previously presented algorithms and running them on a smart charging simulator based on real EV charging data.

%TODO: Stopped reading here

\subsection{Experiment setup}

\subsubsection{Datasets}
We use the Caltech Adaptive charging network dataset~\cite{caltech}. %It contains charging session logs with duration, connection and disconnection time, kWh request, announced parking duration, and slot identifiers. 
The training covers 83 days and contains 1123 sessions across $n=32$ slots.  Each slot has at least 20 training sessions. The test set covers 22 days and 352 sessions over the same number of slots. The dataset is pre-processed to obtain fixed interval data points of $\Delta t=15$ min. Unsatisfiable user requests are rounded down to the largest amount of energy that can be provided in the span of the announced parking time. This means that all requests can be satisfied if the users let their EV connected to the slot for the parking duration that they anticipated. However, this does not always happen in practice. These early disconnections are responsible for \texttt{P-MPC} not being able to reach 100\% satisfaction.

\subsubsection{EVCS settings}
The EVCS parameters are of an existing EVCS station. We set the maximum slot power to 12kW. Each slot can thus deliver an energy amount equal to $\bar{e}=3$kWh in the span of one time step (15 minutes). We further set the subscribed power limit to 8\% of the total nominal capacity of the EVCS, i.e. $\bar{c}= 0.08 \cdot n \cdot \bar{e} = 7.68$kWh per time step. This corresponds to two times the average EVCS load. The electricity prices are set to 0.102€ ~during off-peak hours (i.e. from 00:00 to 06:00, 09:00 to 11:00, 13:00 to 17:00, and 21:00 to 00:00) and to 0.153€ during peak hours. Similarly, the threshold overrun penalty is set to $\xi=14.31$€ for every time step spent above the threshold, as per the ``Tarif Jaune'' rate from France's main electric utility company EDF~\cite{EDF}. Finally, the battery charge efficiency of the connected EVs is $\eta=0.91$

\subsubsection{Implementation specifics}
We run all algorithms on the test dataset with a control horizon $R=40$. \texttt{2S} uses $K'=2$ uncertainty scenarios obtained from a total of $K=20$ samples. We model the optimization problems with Pyomo~\cite{PYOMO} and solve them using HiGHS~\cite{HIGHS}. The classification and regression models described in \ref{sub:uncertainty} are implemented using CatBoost~\cite{prokhorenkova2019catboost}.
\subsubsection{Evaluation Metrics}
The algorithms are compared based on cost and two customer satisfaction criteria:
%The first is their total electricity cost, including contract threshold penalty. The second is the filling rate, i.e., the average percentage of user request provided by session end. The third is the full satisfaction rate, i.e. the percentage of requests fully satisfied by session end.
\begin{itemize}
\item \textit{Electricity Costs}: Total cost, including contract threshold penalties.
\item \textit{Satisfaction Metrics}:
\begin{itemize}
\item \textit{Filling Rate}: Average percentage of user request provided by session end.
\item \textit{Full Satisfaction Rate}: Percentage of requests fully satisfied by session end.
\end{itemize}
\end{itemize}

\begin{figure}[t]
    \centering
    \includegraphics[width=0.8\columnwidth]{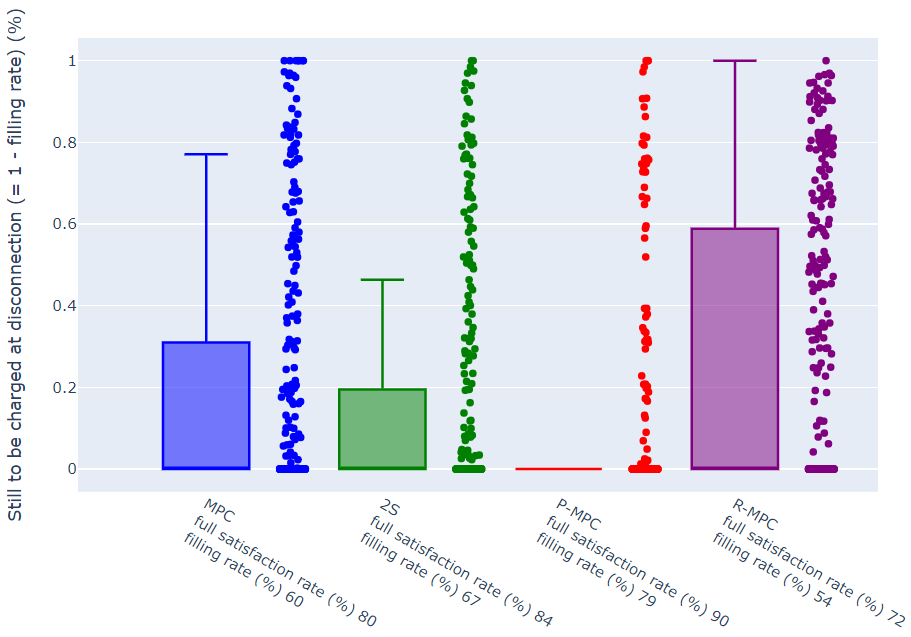}
    \caption{Distribution of the share of the initial energy request that could not be provided by the end of the charging sessions. $\alpha = 5000$}
    \label{fig:distr1}
\end{figure}

\definecolor{PeriwinkleGray}{rgb}{0.749,0.839,0.909}
\definecolor{PastelGreen}{rgb}{0.419,0.901,0.603}
\definecolor{PetiteOrchid}{rgb}{0.85,0.588,0.588}
\begin{table}
\centering
\caption{Simulation results. The results in bold correspond to the version of \texttt{2S} that offers the closest results to \texttt{P-MPC}, while achieving high satisfaction rates}

\label{tab:results}
\resizebox{\linewidth}{!}{%
\begin{tblr}{
  cells = {c},
  cell{1}{1} = {r=2}{},
  cell{1}{2} = {r=2}{},
  cell{1}{3} = {r=2}{},
  cell{3}{5} = {PeriwinkleGray},
  cell{3}{6} = {PeriwinkleGray},
  cell{3}{7} = {PeriwinkleGray},
  cell{3}{8} = {PeriwinkleGray},
  cell{3}{9} = {PeriwinkleGray},
  cell{4}{5} = {PeriwinkleGray},
  cell{4}{6} = {PeriwinkleGray},
  cell{4}{7} = {PeriwinkleGray},
  cell{4}{8} = {PeriwinkleGray},
  cell{4}{9} = {PeriwinkleGray},
  cell{6}{5} = {PastelGreen},
  cell{6}{6} = {PastelGreen},
  cell{6}{7} = {PastelGreen},
  cell{6}{8} = {PastelGreen},
  cell{6}{9} = {PastelGreen},
  cell{7}{5} = {PastelGreen},
  cell{7}{6} = {PastelGreen},
  cell{7}{7} = {PastelGreen},
  cell{7}{8} = {PastelGreen},
  cell{7}{9} = {PastelGreen},
  cell{9}{5} = {PastelGreen},
  cell{9}{6} = {PastelGreen},
  cell{9}{7} = {PastelGreen},
  cell{9}{8} = {PastelGreen},
  cell{9}{9} = {PastelGreen},
  cell{10}{5} = {PastelGreen},
  cell{10}{6} = {PastelGreen},
  cell{10}{7} = {PastelGreen},
  cell{10}{8} = {PastelGreen},
  cell{10}{9} = {PastelGreen},
  cell{12}{5} = {PetiteOrchid},
  cell{12}{6} = {PetiteOrchid},
  cell{12}{7} = {PetiteOrchid},
  cell{12}{8} = {PetiteOrchid},
  cell{12}{9} = {PetiteOrchid},
  cell{13}{5} = {PetiteOrchid},
  cell{13}{6} = {PetiteOrchid},
  cell{13}{7} = {PetiteOrchid},
  cell{13}{8} = {PetiteOrchid},
  cell{13}{9} = {PetiteOrchid},
}
\textbf{$\alpha$} & \textbf{Alg.} &  & \textbf{Electricity} & \textbf{Filling}   & \textbf{Full satisfaction} & \textbf{Electricity cost relative} & \textbf{Filling rate relative}      & \textbf{Full satisfaction relative}           \\
                  &               &                    & \textbf{cost (EUR)}  & \textbf{rate (\%)} & \textbf{rate (\%)}         & \textbf{difference w.r.t P-MPC (\%)} & \textbf{difference w.r.t P-MPC (\%)} & \textbf{difference w.r.t P-MPC (\%)} \\
500               & 2S            &               & 12991                & 69\%               & 60\%                       & 16\%                           & 3\%                             & -1\%                                 \\
500               & MPC           &               & 12392                & 68\%               & 59\%                       & 10\%                           & 1\%                             & -2\%                                 \\
500               & R-MPC         &              & 32942                & 70\%               & 52\%                       & 194\%                          & 4\%                             & -14\%                                \\
1000              & 2S            &              & 17286                & 77\%               & 66\%                       & 3\%                          & -4\%                            & -9\%                                 \\
1000              & MPC           &              & 17156                & 77\%               & 65\%                       & -4\%                           & -5\%                            & -10\%                                \\
1000              & R-MPC         &              & 34519                & 72\%               & 53\%                       & 93\%                           & -10\%                           & -26\%                                \\
5000              & 2S            &              & 43094                & 85\%               & 67\%                     & 3\%                            & -6\%                            & -16\%                                \\
5000              & MPC           &              & 38484                & 81\%               & 61\%                       & -8\%                           & -10\%                           & -23\%                                \\
5000              & R-MPC         &              & 37859                & 73\%               & 55\%                       & -9\%                           & -19\%                           & -31\%                                \\
50000             & 2S            &             & 138893               & 92\%               & 83\%                       & 122\%                          & -1\%                            & -2\%                                 \\
50000             & MPC           &             & 117179               & 89\%               & 78\%                       & 87\%                           & -4\%                            & -7\%                                 \\
50000             & R-MPC         &             & 45470                & 74\%               & 54\%                       & -27\%                          & -21\%                           & -36\%                                
\end{tblr}
}
\end{table}

\begin{figure}[t]
    \centering
\includegraphics[width=0.9\columnwidth]{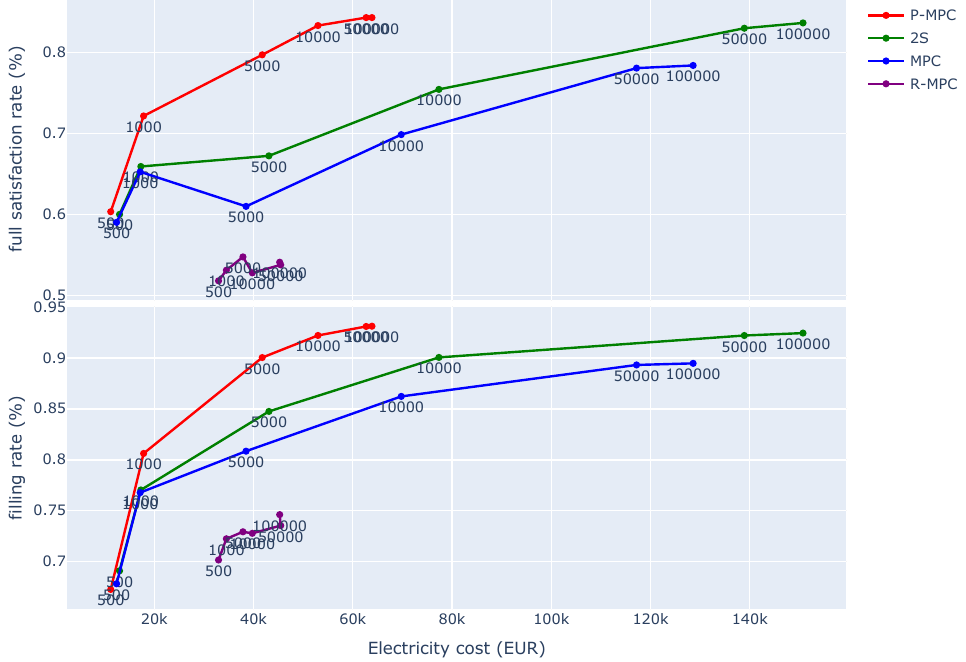}
    \caption{Full satisfaction and filling rates as a function of the electricity costs. text annotations correspond to the values of $\alpha$}
    \label{fig:satisfaction_vs_cost}
\end{figure}
\subsection{Results}

We present the results in Table \ref{tab:results} for varying $\alpha$, representing the weight assigned to customer satisfaction:

% For small $\alpha$ \textcolor{PeriwinkleGray}{$\blacksquare$} (low dissatisfaction penalization), all algorithms except \texttt{R-MPC} exhibit low filling and full satisfaction rates. Notably, \texttt{MPC} and \texttt{2S} incur higher electricity costs compared to \texttt{P-MPC}, attributed to less efficient charging session planning, particularly during off-peak hours.

% As $\alpha$ increases \textcolor{PastelGreen}{$\blacksquare$}, the satisfaction metrics improve across all algorithms. At $\alpha = 5000$, \texttt{2S} closely matches \texttt{P-MPC} in filling rate and full satisfaction rate, outperforming \texttt{R-MPC} by 16\% and 23\% respectively.

% With large $\alpha$ \textcolor{PetiteOrchid}{$\blacksquare$}, \texttt{2S} approaches \texttt{P-MPC}'s satisfaction levels, while \texttt{MPC}'s slightly lag. At $\alpha = 10^5$, \texttt{2S} achieves nearly identical filling and full satisfaction rates as \texttt{P-MPC}, significantly surpassing \texttt{R-MPC} by 24\% and 55\% respectively. However, both \texttt{MPC} and \texttt{2S} incur notably higher electricity costs due to increased dissatisfaction penalization.

For small $\alpha$ \textcolor{PeriwinkleGray}{$\blacksquare$} (low penalization of unsatisfied requests), the filling and full-satisfaction rates are low but identical across all algorithms except \texttt{R-MPC}, which shows lower values. The electricity cost with \texttt{MPC} and \texttt{2S} are larger than those obtained with \texttt{P-MPC}, which could be explained by the fact that these algorithms can not plan the charging sessions as well as \texttt{P-MPC}, and hence can not benefit from the off-peak hours as much. %The electricity cost under the policy of \texttt{R-MPC} is low since a smaller amount of electricity is bought.

As the value of $\alpha$ increases \textcolor{PastelGreen}{$\blacksquare$}, the satisfaction metrics both increase for all algorithms, and our \texttt{2S} and \texttt{MPC} algorithms show similar results to \texttt{P-MPC}. For $\alpha = 5000$, \texttt{2S} is only 3\% more costly but achieves 94\% of the performance of \texttt{P-MPC} in terms of filling rate, and 84\% in terms of full satisfaction rate. These values are respectively 16\% and 23\% higher than the ones obtained with our industry-grade baseline \texttt{R-MPC}.

As the value of $\alpha$ gets very large \textcolor{PetiteOrchid}{$\blacksquare$}, the results obtained by \texttt{2S} are close to the ones obtained with \texttt{P-MPC} in terms of satisfaction, while those obtained with \texttt{MPC} are slightly lower. With $\alpha = 10^5$,  The filling rates of \texttt{2S} and \texttt{P-MPC} are almost identical (92\% for \texttt{2S} against 93\% for \texttt{P-MPC}), and the full satisfaction rates are equal. \texttt{2S} improves the filling rate obtained using \texttt{R-MPC}, our industry-grade baseline by 24\%, and the full satisfaction rate by 55\%. The electricity costs are significantly higher for both \texttt{MPC} and \texttt{2S} (respectively 101\% and 136\% higher than with \texttt{P-MPC}), as customer unsatisfaction is much more penalized than electricity costs. Note that the fact that \texttt{P-MPC} does not reach 100\% of satisfaction rate is because of early disconnections of sessions that could have otherwise been satisfied.

Fig. \ref{fig:satisfaction_vs_cost} illustrates \texttt{2S}'s superior robustness in customer satisfaction compared to \texttt{MPC}, attributed to a broader range of uncertainty scenarios. This enables better anticipation of early session ends, resulting in earlier and faster charging actions. \texttt{2S} consistently outperforms \texttt{MPC} across all satisfaction metrics for equivalent electricity costs, positioning it as the most viable implementation choice. Additionally, \texttt{2S} significantly outperforms \texttt{MPC}, particularly for high $\alpha$ values.

Finally, Fig. \ref{fig:distr1} displays the distribution of unfulfilled energy requests at disconnection time for $\alpha = 5000$. \texttt{R-MPC} exhibits the highest values, reflecting its reliance on user-provided session end times without correction. On average, for $\alpha=5000$, \texttt{P-MPC}, \texttt{2S}, \texttt{MPC}, and \texttt{R-MPC} each took 164ms, 1563ms, 556ms, and 193ms per time step on an Intel i7 1185G7 laptop.

\section{Conclusion}
In this paper, we presented an Electric Vehicle Charging Station model that considers real-world constraints, including slot power limitations, penalties for contract threshold overruns, and early disconnections of EVs. We proposed a mathematical formulation for the optimal control problem under uncertainty and implemented two Multi-Stage Stochastic Programming approaches: Model Predictive Control (\texttt{MPC}), and Two-Stage Stochastic Programming  (\texttt{2S}). The problem complexity arises from uncertainties in charging session start and end times, as well as the energy requested by users. This necessitates predicting future charging session times and their energy demand. Early disconnections pose an additional challenge, as customers may disconnect their EVs before the initially announced time, further complicating charge planning. To address this, we introduced a user behavior model based on a sojourn time-dependent stochastic process that considers the duration of the charging session, allowing for cost reduction while maintaining customer satisfaction.

Our algorithms were compared to two baselines: a Model Predictive Control algorithm with a perfect forecast of the future uncertainties (\texttt{P-MPC}), representing a theoretical optimal controller, and an industry-grade Model Predictive Control algorithm (\texttt{R-MPC}), mentioned in~\cite{tristanspaper}. Simulations were conducted using a real-world dataset spanning 22 days. The results demonstrated the advantages of our methods, which incorporate user behavior modeling and leverage user-provided information. Particularly, \texttt{2S} proves robust against early disconnections by considering a greater number of uncertainty scenarios for optimization. The version of \texttt{2S} that prioritizes user satisfaction over electricity cost performs comparably to the optimal baseline \texttt{P-MPC}, outperforming our previous algorithm from~\cite{tristanspaper}, \texttt{R-MPC}, by 24\% and 55\% in the two user satisfaction metrics, respectively. The version of \texttt{2S} that strikes the best balance between cost and user satisfaction shows a 3\% relative cost difference (slightly higher) compared to the optimal baseline \texttt{P-MPC}, while achieving 94\% and 84\% of the user satisfaction performance of \texttt{P-MPC} in the two metrics used. \newline\newline

\bibliographystyle{IEEEtran}
\bibliography{references}  

\end{document}